\documentclass[12pt,a4paper]{article}
\usepackage[utf8]{inputenc}
\usepackage[english]{babel}
\usepackage{amsmath}
\usepackage{amssymb}
\usepackage{amsthm} 
\newtheorem{theorem}{Theorem}

\newtheorem{remark}{Remark}

\usepackage[svgnames]{xcolor}
\def\ep{\varepsilon}

\usepackage{lineno}

\title{$L_2$-small ball asymptotics \\
for some demeaned Gaussian processes}
\date{\today}
\author{
Alexander Nazarov\footnote{St.Petersburg Dept of Steklov Institute, Fontanka 27, St.Petersburg, 191023, Russia, 
    and St.Petersburg State University, Universitetskii pr. 28, St.Petersburg, 198504, Russia. E-mail: al.il.nazarov@gmail.com. Supported by Russian Scientific Foundation, Grant 21-11-00047.} \ and \setcounter{footnote}{6}
Yulia Petrova\footnote{Pontifícia Universidade Católica do Rio de Janeiro (PUC-Rio), R. Marquês de São Vicente, 124, Gávea, Rio de Janeiro, 22451-040, Brazil. E-mail: yu.pe.petrova@yandex.ru.}}

\begin{document}
\maketitle
\section{Introduction. Main result}
The theory of small ball probabilities (also called small deviation probabilities) is extensively studied in recent decades (see the surveys~\cite{Lf}, \cite{LS}, \cite{F}; for the extensive up-to-date bibliography see~\cite{L-table}). The most explored case is that of Gaussian processes in Hilbertian norm, see the recent survey~\cite{L2-survey}. Namely, for
a zero mean-value Gaussian process $X(t)$, $0\leq t\leq 1$, we are interested in the sharp $L_2$-small ball asymptotics, that is the relation
$$
\mathbb{P}\Big\{\|X\|^2:= \int\limits_0^1 X^2(t)\,dt\le\varepsilon^2\Big\}\sim f(\varepsilon),\quad \varepsilon\to 0.
$$
By the Karhunen--Lo{\`e}ve expansion (see, e.g.,~\cite[Chapter 2]{Lif2012-book}) we have the following distributional equality:
\begin{align}
    \label{eq:KL}
    \|X\|^2\stackrel{d}{=}\sum\limits_{k=1}^{\infty}\mu_k\xi_k^2,
\end{align}
where $\mu_k$, $k\in\mathbb{N}$, stand for the positive eigenvalues of the integral equation
\begin{align}
    \label{eq:X-eigen}
\int\limits_0^1 G(t,s)u(s)\,ds=\mu u(t),
\end{align}
counted in non-increasing order with their multiplicities, while $\xi_k$ are standard i.i.d. normal random variables.
The kernel $G$ in \eqref{eq:X-eigen} is the covariance function of the process $X$: 
$$
G(t,s)\equiv G_X(t,s):=\mathbb{E}X(t)X(s). 
$$

Formula \eqref{eq:KL} shows that the $L_2$-norm distribution of a Gaussian process $X$ depends only on the eigenvalues of \eqref{eq:X-eigen}. Moreover, the classical Wenbo Li comparison principle \cite{Li} claims that to obtain the sharp (at least up to a multiplicative constant) small ball asymptotics for $X$ we need sufficiently sharp (usually two-terms) eigenvalues asymptotics of the corresponding integral equation. Most of the results in this direction were made for the special important class of Gaussian processes introduced in \cite{NN}.

Recall that the \textit{\textbf{Green Gaussian process}} is a Gaussian process $X$ on the interval (say, $[0,1]$) such that its covariance function $G_X$ is the (generalized) Green function of 
the boundary value problem
\begin{align}
\label{eq:L}
{\cal L}u:=
(-1)^{\ell}\left(p_{\ell}u^{({\ell})}\right)^{({\ell})}+
\left(p_{\ell-1}u^{({\ell}-1)}\right)^{({\ell}-1)}+\dots+p_0u=\lambda u
\end{align}
(here $p_j$, $j=0,\dots,\ell$ are functions on $[0,1]$, and
$p_{\ell}(t)>0$) with proper boundary conditions.\footnote{Definitions of the Green function and generalized Green function can be found, e.g., in~\cite[Ch. 2, \S1]{Sm}, see also~\cite{Naz09a}. For simplicity we assume that $p_j\in W^j_\infty[0,1]$, $j=0,\dots,\ell$.} This class of processes is very important as it includes many classical processes -- the Wiener process, the Brownian bridge, the Ornstein-Uhlenbeck process, their (multiply) integrated counterparts etc. It is easy to see that for a Green Gaussian process the (non-zero) eigenvalues $\mu_k$ in \eqref{eq:KL} are related to the (non-zero) eigenvalues $\lambda_k$ of the boundary value problem \eqref{eq:L} by relation $\mu_k=\lambda_k^{-1}$. So, to obtain rather good asymptotics of eigenvalues $\mu_k$ one can use the powerful methods of spectral
theory of ordinary differential operators, see, e.g., \cite[Sec. 4]{Nm} and \cite[Chap. XIX]{DS}.
\medskip

In this paper we deal with the class of centered (\textit{\textbf{demeaned}}) processes. The operation of centering is classical and dates back at least to Watson \cite{W}, who used it for testing nonparametric hypotheses on the circle. 
For a random process $X(t)$ on $[0,1]$ we define the process 
\begin{align}
\label{eq:demeaned}
\overline{X}(t)=X(t)-\int\limits_0^1 X(s)\,ds.
\end{align}
The following problem arises naturally in this relation: when the operation of centering transforms a Green Gaussian process to a Green one. For some other natural operations on Gaussian processes this problem was investigated in a number of papers, see~\cite{NN}, \cite{Naz09a}, \cite{NP}, \cite{NN21}, and \cite[Sec. 4.1]{L2-survey}.

It is easy that if $G(t,s)$ is the covariance function of $X$ then the covariance function of $\overline{X}$ is
\begin{align}
\label{eq:G-demeaned}
\overline{G}(t,s)
:=G(t,s)
-\int\limits_0^1G(t,y)\,dy
-
\int\limits_0^1G(x,s)\,dx+
\int\limits_0^1\int\limits_0^1G(x,y)\,dydx.
\end{align}
The covariance operator $\overline{\cal G}$ with demeaned kernel~\eqref{eq:G-demeaned} clearly has a zero eigenvalue with a constant eigenfunction, so it cannot be the Green function of a boundary value problem but can be a generalized Green function.
\medskip

Our main result is the following theorem. Its particular case was proved in~\cite[Theorem 3.1, part 2]{Naz09a}.
\begin{theorem}
\label{thm:Green-to-Green-demeaned}
Let $X$ be a Green Gaussian process on $[0,1]$, and let the corresponding differential operator ${\cal L}$ have no zero order term ($p_0\equiv0$).   
Then the demeaned process \eqref{eq:demeaned} is also the Green Gaussian process. Namely, the covariance function~\eqref{eq:G-demeaned} is the generalized Green function for the boundary value problem $\overline{\cal L}u:={\cal L}u=\lambda u$ with boundary conditions specified below.
\end{theorem}

\begin{remark}
    The case $p_0\not\equiv0$ is more complicated. In this case the demeaned process~\eqref{eq:demeaned} is in general not a Green Gaussian process. 
\end{remark}

\begin{proof} 
If the covariance operator  ${\cal G}$ of the original process $X$ has zero eigenvalue with a constant eigenfunction then easily $\overline{G}\equiv G$ and $\overline{X}=X$. So, we assume that this is not the case. For simplicity only, we assume also that ${\cal G}$ has no non-constant functions corresponding to the zero eigenvalue. In this case $G(t,s)$ is the conventional Green function, that is ${\cal L}G(t,s)=\delta(t-s)$. Since $p_0\equiv0$, straightforward calculation yields
\begin{align}
\label{eq:delta-1}
{\cal L}\,\overline{G}(t,s)=\delta(t-s)-1.
\end{align}
So, to complete the proof we need only to define boundary conditions for $\overline{\cal L}$.

The boundary conditions for ${\cal L}$ can be written as follows: 
\begin{align}
\label{eq:bc-non-separated}
U_{\nu}(u):= U_{\nu 0}(u)+U_{\nu 1}(u)=0,\qquad \nu=1,\dots,2\ell, 
\end{align}
where\footnote{Usually the boundary conditions are supposed to be Birkhoff-normalized, see, e.g., \cite[Section 2, 4.7]{Nm} and~\cite{Naz09a}. In particular, they are written with separated highest order terms. For our purposes it is more convenient to separate the zero order terms. Thus, we do not suppose a priori that the boundary conditions \eqref{eq:bc-non-separated} are normalized.}
\begin{align*}
U_{\nu 0}(u):=\alpha_{\nu 0}u(0)+\sum\limits_{j=1}^{2\ell-1}
\alpha_{\nu j}u^{(j)}(0),
\qquad
U_{\nu 1}(u):=\gamma_{\nu0}u(1)+\sum\limits_{j=1}^{2\ell-1}
\gamma_{\nu j}u^{(j)}(1).
\end{align*}

By linear transformation we can rewrite the boundary conditions~\eqref{eq:bc-non-separated} in an equivalent form such that not more than two forms $U_\nu(u)$ contain the zero order terms. Without loss of generality they correspond to $\nu=1$ and $\nu=2$, that is, $\alpha_{\nu0}=\gamma_{\nu0}=0$ for $\nu>2$. We should distinguish three cases.
\medskip

{\bf 1}. If non of the forms $U_{1}(u)$ and $U_{2}(u)$ contains the zero order terms, then the boundary value problem~\eqref{eq:L},~\eqref{eq:bc-non-separated} has a constant eigenfunction corresponding to the zero eigenvalue, that is excluded from the very beginning.\medskip

{\bf 2}. If both forms $U_{1}(u)$ and $U_{2}(u)$ contain the zero order terms, they can be rewritten as follows:\footnote{If this is impossible then we can in fact eliminate the zero order terms in the form $U_2(u)$ and arrive at the next case.}
\begin{equation}
\label{eq:U1-U2}
\gathered
U_{1}(u):=u(0)+
\sum\limits_{j=1}^{2\ell-1}
\alpha_{1j}u^{(j)}(0)+\sum\limits_{j=1}^{2\ell-1}
\gamma_{1j}u^{(j)}(1)=0,
\\
U_{2}(u):=u(1)+
\sum\limits_{j=1}^{2\ell-1}
\alpha_{2j}u^{(j)}(0)+\sum\limits_{j=1}^{2\ell-1}
\gamma_{2 j}u^{(j)}(1)=0.
\endgathered
\end{equation}
It is easy to see that the function $\overline{G}(t,s)$ satisfies all the same boundary conditions as $G(t,s)$ except for \eqref{eq:U1-U2}. However, since the subtraction of a constant keeps the difference $u(0)-u(1)$, $\overline{G}(t,s)$ satisfies the boundary condition
\begin{align*}
U_2(u)-U_1(u)=0.
\end{align*}

Finally, we notice that the operator ${\cal L}$ can be rewritten as ${\cal L}u=(Lu'){}'$, where $L$ is a $(2\ell-2)$-order operator. Therefore, formula~\eqref{eq:delta-1} implies
\begin{align*}
0=\int\limits_0^1 {\cal L}\overline{G}(y,s)\,dy
=
\int\limits_0^1 (L \overline{G}{}^{\,\prime}_y(y,s)){}^{\prime}_y\,dy
=
(L\overline{G}{}^{\,\prime}_t)(t,s)\Big|_{t=0}^{t=1},
\end{align*}
and the last boundary condition is
\begin{align}
\label{eq:bc-demeaned}
(Lu')(0)-(Lu')(1)=0.
\end{align}

{\bf 3}. Let exactly one of the forms, say $U_{1}(u)$, contain the zero order terms:
\begin{equation}
\label{eq:U1}
U_{1}(u):=\alpha_{10} u(0)+\gamma_{10}u(1)+
\sum\limits_{j=1}^{2\ell-1}
\alpha_{1j}u^{(j)}(0)+\sum\limits_{j=1}^{2\ell-1}
\gamma_{1j}u^{(j)}(1)=0.
\end{equation}
If $\alpha_{10}+\gamma_{10}=0$ then again the boundary value problem~\eqref{eq:L},~\eqref{eq:bc-non-separated} has a constant eigenfunction corresponding to the zero eigenvalue, that is excluded. 

Otherwise, $\overline{G}(t,s)$ satisfies all the same boundary conditions as $G(t,s)$ except for \eqref{eq:U1}. Then, as in the case {\bf 2}, the last boundary condition is~\eqref{eq:bc-demeaned}.
\end{proof}

In the next Section we consider two applications of Theorem \ref{thm:Green-to-Green-demeaned}. Also we consider two examples where this theorem is not applicable but nevertheless sharp $L_2$-small ball asymptotics for corresponding demeaned processes can be obtained.

\section{Examples}

{\bf Example 1}. 
Consider a zero-mean process $X_\alpha(t)$, which covariance function $G_{X_\alpha}$ is a Green function of the following boundary value problem:
\begin{align*}
{\cal L}u\equiv-u''=\lambda u,\qquad u(0)+\alpha u(1)=0, \qquad \alpha u'(0)+u'(1)=0.
\end{align*}

If $\alpha=-1$ (periodic boundary conditions) then $X_\alpha$ has zero eigenvalue with a constant eigenfunction (in fact, $X_\alpha$ coincides in distribution with the centered Brownian bridge $\overline{B}$). So, we assume that $\alpha\ne -1$.

For $\alpha=0$, $X_\alpha$ is the standard Wiener process $W(t)$; the case $\alpha=1$ (anti-periodic boundary conditions) corresponds to the generalized Slepian process $S_{\frac12}(t)$ (see, e.g., \cite{GL}, \cite{Naz09a}). According to Theorem~\ref{thm:Green-to-Green-demeaned}, for the demeaned process $\overline{X_\alpha}$ the corresponding boundary value problem reads as follows:
\begin{align*}
{\cal L}u=-u''=\lambda u,\qquad u'(0)=u'(1)=0,
\end{align*}
It is worth to note that there is no dependence on $\alpha$. Therefore, we have equality in distribution
$$
\|\overline{X_\alpha}\|\stackrel{d}{=}\|\overline{W}\|.
$$
In particular, 
$$
    \mathbb{P}\Big\{\|\overline{X_\alpha}\|\leq\ep\Big\}
        =
    \mathbb{P}\Big\{\|\overline{W}\|\leq \ep\Big\}\sim
    \Big(\frac 8{\pi}\Big)^{\frac 12}\exp\Bigl(-\frac18\,\ep^{-2}\Bigr)
$$
(the latter relation was obtained in~\cite{BNO}).
\medskip

{\bf Example 2}. 
Now we consider integrated Ornstein-Uhlenbeck process
$$
(U_\beta)_1(t)=\int\limits_0^t U_\beta(\tau)\,d\tau,
$$
where $U_\beta$, $\beta>0$, is the standard Ornstein-Uhlenbeck (OU) process with covariance function 
$$
G_{U_\beta}(t,s)=\frac 1{2\beta}\exp(-\beta|t-s|).
$$
It is well known that $U_\beta$ is a Green Gaussian process corresponding to the following boundary value problem:
\begin{equation}
\label{eq:OU}
Lu\equiv -u''+\beta^2 u=\lambda u,  \quad u'(0)-\beta u(0)=0, \quad u'(1)+\beta u(1)=0.
\end{equation}

By \cite[Theorem 2.1]{NN}, the covariance of integrated OU process is the Green function of the boundary value problem
\begin{equation}
\label{eq:OU-int}
\begin{aligned}
{\cal L}u\equiv u^{IV}-\beta^2 u''=\lambda u, & \qquad u(0)=0, \quad u'''(1)-\beta^2u'(1)=0,\\
u''(0)-\beta u'(0)=0, & \qquad u''(1)+\beta u'(1)=0.
\end{aligned}
\end{equation}
By Theorem 1, the boundary value problem corresponding to the demeaned process $\overline{(U_\beta)_1}$ reads as follows:
\begin{align}
\nonumber
{\cal L}u\equiv u^{IV}-\beta^2 u''=\lambda u, 
    \\ 
\label{eq:BVP-demeaned-OU}
    u'''(0)-\beta^2u'(0)=0, & \qquad u'''(1)-\beta^2u'(1)=0,
    \\
    \nonumber
u''(0)-\beta u'(0)=0, & \qquad u''(1)+\beta u'(1)=0.
\end{align}

By Theorem 7.1 in \cite{NN}, the eigenvalues of this problem have the following two-term asymptotics:
$$
\lambda_k=\pi^4\big(k-\frac 32+O(k^{-1})\big)^4, 
\quad
k\to\infty.
$$

Since one of the eigenvalues of the problem \eqref{eq:OU-int} equals zero (it corresponds to the constant eigenfunction), we should shift the numeration and obtain the following asymptotics of $\mu_k$:
$$
\mu_k=\lambda_{k+1}^{-1}=\frac 1{\pi^4}\big(k-\frac 12+O(k^{-1})\big)^{-4}, 
\quad
k\to\infty.
$$
Using the Wenbo Li comparison principle and Theorem 6.2 in \cite{NN}, we obtain, as $\ep\to0$,
\begin{align}
\label{eq:intOU}    
P\{||\overline{(U_\beta)_1}||
\leq\ep\} \sim C_{\rm dist}\cdot \frac {8}{(3\pi)^{\frac 12}}
\cdot
\ep^{\frac 13} \exp \big(- \frac 38\, \ep^{-\frac 23} \big),
\end{align}
where the distortion constant is defined by
$$
C_{\rm dist}=\prod_{k=1}^{\infty}\,
\frac{\lambda_{k+1}^{\frac 12}}{{\pi^2}(k-\frac 12)^2 }.
$$

To compute $C_{\rm dist}$, we use the algorithm developed in \cite[Section 2]{Na}\footnote{In fact, in the case under consideration some arguments could be simplified. However, we provide them in the general form.}.
Let $\pm\zeta_j$, $j=0,1$, be the complex roots of the equation
\begin{align}
\label{eq:eigen}    
\zeta^{4}+\beta^2\zeta^{2}=\lambda.
\end{align}
Denote by $\zeta_0$ the root of \eqref{eq:eigen} which is positive for positive $\lambda$ (obviously, when $\lambda>0$, \eqref{eq:eigen} has just one pair of real roots). 
Then we can consider $\zeta_1$ as the function of $\zeta=\zeta_0$.

Next, when all the roots of \eqref{eq:eigen} are different the functions 
$e^{\pm i\zeta_jt}$, $j=0,1$, form the fundamental system of 
the equation ${\cal L}u=\lambda u$. 
Substituting them to the boundary conditions we obtain that if $\lambda$ is an eigenvalue of \eqref{eq:BVP-demeaned-OU} then the roots
of \eqref{eq:eigen} are the zeros of the determinant
\begin{align}
\label{eq:det}
{\scriptsize
{\mathfrak F}(\zeta)
=
\det
\begin{bmatrix}
-i\zeta_0^3 -\beta^2 i\zeta_0&i\zeta_0^3+\beta^2 i \zeta_0&-i\zeta_1^3 -\beta^2 i\zeta_1&i\zeta_1^3+\beta^2 i \zeta_1
\\
-\zeta_0^2-\beta i\zeta_0&-\zeta_0^2+\beta i \zeta_0& -\zeta_1^2-\beta i\zeta_1&-\zeta_1^2+\beta i \zeta_1
\\
(-i\zeta_0^3 -\beta^2 i\zeta_0)e^{i\zeta_0}&(i\zeta_0^3+\beta^2 i \zeta_0)e^{-i\zeta_0}&(-i\zeta_1^3 -\beta^2 i\zeta_1)e^{i\zeta_1}&(i\zeta_1^3+\beta^2 i \zeta_1)e^{-i\zeta_1}
 \\
(-\zeta_0^2+\beta i\zeta_0)e^{i\zeta_0}&(-\zeta_0^2-\beta i \zeta_0)e^{-i\zeta_0}& (-\zeta_1^2+\beta i\zeta_1)e^{i\zeta_1}&(-\zeta_1^2-\beta i \zeta_1)e^{-i\zeta_1}
\end{bmatrix}.
}
\end{align}

When two roots of \eqref{eq:eigen} coincide, ${\mathfrak F}(\zeta)$ also vanishes. To eliminate these irrelevant zeros we consider the function
$$
{\mathfrak F}_1(\zeta)= \frac {{\mathfrak F}(\zeta)}
{{\mathfrak V}(\zeta_0,-\zeta_0,\zeta_1,-\zeta_1)},
$$
where ${\mathfrak V}$ stands for the Vandermonde determinant.
Now, $\lambda$ is an eigenvalue of \eqref{eq:BVP-demeaned-OU} iff the roots of \eqref{eq:eigen}
are the zeros of ${\mathfrak F}_1(\zeta)$.

Since we should exclude the eigenvalue $\lambda=0$, we introduce the function 
$$
{\mathfrak F}_2(\zeta)= \frac {{\mathfrak F}_1(\zeta)}
{\zeta^4+\beta^2\zeta^2}.
$$

Note that as $|\zeta|\to\infty$, we have
\begin{align}
    \label{eq:zeta-1}
    \zeta_1=i\zeta+O(|\zeta|^{-1}).
\end{align}

\noindent Therefore, as $|\zeta|\to\infty$,
$$
|{\mathfrak V}(\zeta_0,-\zeta_0,\zeta_1,-\zeta_1)|\sim
|\zeta|^6\cdot|{\mathfrak V}(1,-1,i,-i)|=
16|\zeta|^6.
$$ 

Next, due to \eqref{eq:zeta-1} we have, see \cite[\S 4, Theorem 2]{Nm}, as $|\zeta|\to\infty$ 
and $|\arg(\zeta)|\le \pi/4$, 
$$
{\mathfrak F}(\zeta)=\zeta^{10}\exp(\zeta)
\cdot
\Bigl(\Phi(\zeta)+O(|\zeta|^{-1})\Bigr),
$$
where
$$
\Phi(\zeta)=\det
\begin{bmatrix}
-i&i&-1&0\\ 
-1&-1&1&0\\
-ie^{i\zeta}&ie^{-i\zeta}& 0 & 1 \\ 
-e^{i\zeta}&-ie^{-i\zeta}& 0 & 1 
\end{bmatrix} = 4i\cos(\zeta).
$$
Similar asymptotics holds in other quarters of the complex plane. Therefore, we have, as $|\zeta|\to\infty$, 
$$
|{\mathfrak F}_2(\zeta)|=\frac 12\,\bigl(|\cosh(\zeta)|+O(|\zeta|^{-1})\bigr)
\cdot
\bigl(|\cos(\zeta)|+O(|\zeta|^{-1})\bigr).
$$

Now we consider the auxiliary function 
$$
{\Psi}(\zeta)= \frac 12 \cosh(\zeta)\cos(\zeta).
$$ 
It is easy to see that if we take 
$|\zeta|=\pi N$, then as $N\to \infty$,
\begin{align}
\label{eq:F-psi}
    \frac {|{\mathfrak F}_2(\zeta)|}{|{\Psi}(\zeta)|}
\rightrightarrows 1.
\end{align}
Therefore, by the Rouch\'{e} Theorem \cite[\S3.4.2]{T}, functions ${\mathfrak F}_2$ and $\Psi$ have equal number of zeros in the disk $|\zeta|<\pi N$ for large $N$. We notice that the zeros of $\Psi$ are just $\pm \pi (k-\frac 12)$ and $\pm \pi i (k-\frac 12)$, $k\in\mathbb N$, and denote by $\pm\zeta_{0k}$, $\pm\zeta_{1k}$, $k\in\mathbb N$, the zeros of ${\mathfrak F}_2$. Since $|\zeta_{0k}|^2|\zeta_{1k}|^2=\lambda_{k+1}$ by the Vieta theorem, we have
$$
C^2_{\rm dist}=\prod_{k=1}^{\infty}\,
\frac{|\zeta_{0k}|^2|\zeta_{1k}|^2}{{\pi^4}(k-\frac 12)^4}.
$$
Now we use the Jensen Theorem \cite[\S3.6.1]{T}, push $N\to\infty$ and take into account \eqref{eq:F-psi}. This gives
\begin{align*}
C^2_{\rm dist}=\frac {|{\mathfrak F}_2(0)|}{|\Psi(0)|}=2|{\mathfrak F}_2(0)|.
\end{align*}

Since $\zeta_1\to \beta i$ as $\zeta\to0$, we have
$$
|{\mathfrak V}(\zeta_0,-\zeta_0,\zeta_1,-\zeta_1)|\sim 4\beta^5|\zeta|, \qquad \zeta\to0.
$$
So,
\begin{align}
\label{eq:prod}
C^2_{\rm dist}=\lim\limits_{\zeta\to0}\frac {|{\mathfrak F}(\zeta)|}{2\beta^7|\zeta|^3}.
\end{align}
Adding the second column in \eqref{eq:det} to the first one we obtain
$$
{\mathfrak F}(\zeta)\sim\det
\begin{bmatrix}
0 & \beta^2 i \zeta & 0 & 0
\\
0 & \beta i \zeta & 2\beta^2  & 0
\\
2\beta^2 \zeta \sin(\zeta) & \beta^2 i \zeta e^{-i\zeta} & 0 & 0
 \\
-2\beta \zeta \sin(\zeta) & -\beta i \zeta e^{-i\zeta} & 0 & 2\beta^2e^\beta
\end{bmatrix}
 = 8i\beta^8e^\beta \zeta^3, \quad \zeta\to0.
$$
Taking into account \eqref{eq:intOU} and \eqref{eq:prod}, we obtain finally
$$
P\{||\overline{(U_\beta)_1}||\leq\ep\} \sim 16 \Big(\frac {\beta e^\beta}{3\pi}\Big)^{\frac 12}
\cdot
\ep^{\frac 13} \exp \big(- \frac 38\, \ep^{-\frac 23} \big), \qquad \ep\to0.
$$

\begin{remark}
    In a similar way one can obtain $L_2$-small ball asymptotics for the demeand $m$-times integrated OU process $\overline{(U_\beta)_m}$ for any $m\in\mathbb{N}$.
\end{remark}

{\bf Example 3}. 
Here we consider the conventional OU process $U_\beta$. Notice that the operator $L$ in \eqref{eq:OU} has a non-zero zero order term. So, Theorem \ref{thm:Green-to-Green-demeaned} is not applicable.

The spectrum of the demeaned OU process $\overline{U_\beta}$ was studied in the paper \cite{Ai}. However, the results of \cite{Ai} contain some algebraic errors.

Denote by $\overline{G}$ the Green function of the process $\overline{U_\beta}$. Then the straightforward calculation yields
\begin{align*}
L\overline{G}(t,s)=\delta(t-s)-1- \beta^2\int\limits_0^1 G_{U_\beta}(x,s)\,dx + \beta^2\int\limits_0^1\int\limits_0^1G_{U_\beta}(x,y)\,dydx,
\end{align*}
and thus
\begin{align}
\label{eq:delta'}
(L\overline{G}(t,s))'=\delta'(t-s).
\end{align}

Let $v$ be a non-constant eigenfunction of the covariance operator $\overline{\cal G}$: 
\begin{align}
\label{eq:G-eigen}
\int\limits_0^1\overline{G}(t,s)v(s)\,ds=\mu v(t).
\end{align}
We multiply~\eqref{eq:delta'} by $v(s)$ and integrate over $[0,1]$. This gives
\begin{align}
\label{eq:L'-eigen}
(Lv)'\equiv -v'''+\beta^2v'=\lambda v',\qquad \lambda=\mu^{-1}.
\end{align}
A general solution of this equation is
\begin{align}
\label{eq:sol}
v(t)=C_1+C_2\cos(\zeta t)+C_3\sin(\zeta t), \qquad \lambda=\zeta^2+\beta^2.
\end{align}
Direct substitution of this ansatz into \eqref{eq:G-eigen} gives, after tedious calculations, the homogeneous system of linear algebraic equations for $C_1$--$C_3$. The determinant of this system should vanish, and we obtain
\begin{multline*}
F_\beta(\zeta)\equiv \frac 1{\zeta^4}\,\Big(2\beta(\zeta^2+\beta^2)-2\beta\cos(\zeta)(\zeta^2+\beta\zeta^2+\beta^2)
\\
+\zeta\sin(\zeta)(2\zeta^2+\beta\zeta^2+2\beta^2-\beta^3)\Big)=0.
\end{multline*}

Now we consider the auxiliary function 
$$
{\psi}(\zeta)=(2+\beta)\,\frac {\sin(\zeta)}{\zeta}.
$$ 
It is easy to see that if we take $|\zeta|=\pi (N-1/2)$, 
then as $N\to \infty$,
\begin{align}
\label{eq:F-psi-1}
    \frac {|F_\beta(\zeta)|}{|{\psi}(\zeta)|}
\rightrightarrows 1.
\end{align}
By the Rouch\'{e} Theorem, functions $F_\beta$ and $\psi$ have equal number of zeros in the disk $|\zeta|<\pi (N-1/2)$ for large $N$. We notice that the zeros of $\psi$ are just $\pm \pi k$, $k\in\mathbb N$, and denote by $\pm\zeta_k$, $k\in\mathbb N$, the zeros of $F_\beta$. Then we use the Jensen Theorem, take into account \eqref{eq:F-psi-1} and push $N\to\infty$. This gives
\begin{align}
\label{eq:prod1}
\prod_{k=1}^{\infty}\frac 
{\zeta_k^2} {(\pi k)^2}
=\frac {F_\beta(0)}{{\psi}(0)}.
\end{align}

Now we recall that the non-zero eigenvalues of \eqref{eq:G-eigen} are related to $\zeta_k$ by the formula $\mu_k=\frac 1{\zeta_k^2+\beta^2}$.
In view of \eqref{eq:prod1} we have 
$$
\prod_{k=1}^{\infty}\, \frac {(\pi k)^{-2}}{\mu_k} =
\prod_{k=1}^{\infty}\, \frac {\zeta_k^2+\beta^2}{(\pi k)^2}=\frac {F_\beta(0)}{{\psi}(0)}
\cdot \prod_{k=1}^{\infty}\Big(1+\frac {\beta^2}
{\zeta_k^2}\Big) \stackrel{(*)}{=} \frac {F_\beta(\beta i)}{\psi(0)}=\frac{2e^\beta}{2+\beta}
$$
(the equality $(*)$ follows from the Hadamard Theorem on canonical product \cite[\S8.2.4]{T}).

By the Wenbo Li comparison principle we have, as $\ep\to0$,
$$
 \mathbb{P}\Big\{\sum_{k=1}^{\infty} \mu_k \xi_k^2\leq\ep^2\Big\}
   \sim
 \Big(\frac{2e^\beta}{2+\beta}\Big)^{\frac 12}\,
 \mathbb{P}\Big\{\sum_{k=1}^{\infty} (\pi k)^{-2} \xi_k^2\leq\ep^2\Big\},
$$
or
$$
    \mathbb{P}\Big\{\|\overline{U_\beta} \|\leq\ep\Big\}
        \sim
    \Big(\frac{2e^\beta}{2+\beta}\Big)^{\frac 12}\,
    \mathbb{P}\Big\{\|B\|\leq \ep\Big\}\sim
    \Big(\frac{16e^\beta}{\pi(2+\beta)}\Big)^{\frac 12}\,\exp\Bigl(-\frac18\,\ep^{-2}\Bigr)
$$
(here $B$ is the Brownian bridge, and the latter relation is classical, see, e.g.,~\cite{AD}).
\medskip

{\bf Example 4}. 
Finally, we consider the OU process starting at zero $\stackrel{\circ}U_\beta$ with covariance function 
$$
G_{\stackrel{\circ}U_\beta}(t,s)=\frac 1{2\beta}\big(\exp(-\beta|t-s|)-\exp(-\beta(t+s))\big).
$$ 
It is easy to see that $\stackrel{\circ}U_\beta$ is a Green Gaussian process corresponding to the following boundary value problem:
\begin{equation}
\label{eq:OU-0}
Lu\equiv -u''+\beta^2 u=\lambda u,  \quad u(0)=0, \quad u'(1)+\beta u(1)=0.
\end{equation}
As in the previous example, the operator $L$ in \eqref{eq:OU-0} has a non-zero zero order term, and Theorem \ref{thm:Green-to-Green-demeaned} is not applicable.

The straightforward calculation shows that the Green function of the process $\overline{\stackrel{\circ}U_\beta}$ satisfies the equation \eqref{eq:delta'}. Therefore, any non-constant eigenfunction $v$ of the corresponding covariance operator is a solution of the equation \eqref{eq:L'-eigen}. Direct substitution of \eqref{eq:sol} into the integral equation gives the homogeneous system of linear algebraic equations for $C_1$--$C_3$. Since this system should degenerate, we obtain
\begin{multline*}
F^\circ_\beta(\zeta)\equiv \frac 1{\zeta^4}\,\Big(2\beta(\zeta^2+\beta^2)-\beta\cos(\zeta)(2\zeta^2+\beta\zeta^2+2\beta^2)
\\
+\zeta\sin(\zeta)(\zeta^2+\beta^2-\beta^3)\Big)=0.
\end{multline*}

Now we consider the auxiliary function 
$$
\psi^\circ(\zeta)=\frac {\sin(\zeta)}{\zeta}
$$ 
and, arguing as in previous example, obtain
$$
\prod_{k=1}^{\infty}\, \frac {(\pi k)^{-2}}{\mu_k} = \frac {F^\circ_\beta(\beta i)}{\psi^\circ(0)}=e^\beta.
$$
By the Wenbo Li comparison principle we have, as $\ep\to0$,
$$
\mathbb{P}\Big\{\|\overline{\stackrel{\circ}U_\beta} \|\leq\ep\Big\}
        \sim
    e^{\frac \beta 2}\,
    \mathbb{P}\Big\{\|B\|\leq \ep\Big\}\sim
    \Big(\frac{8e^\beta}{\pi}\Big)^{\frac 12}\,\exp\Bigl(-\frac18\,\ep^{-2}\Bigr).
$$

\begin{remark}
    For $\beta=0$, we have the following distributional equality:
$$
\|\overline{U_0}\|\stackrel{d}{=}\|\overline{\stackrel{\circ}U_0}\|\stackrel{d}{=}\|\overline{W}\|.
$$
\end{remark}


\small
\begin{thebibliography}{GHLT2}

\bibitem{Ai}
X.~Ai, {\em A note on Karhunen–Loeve expansions for the demeaned stationary Ornstein–Uhlenbeck process}, Stat. Probab. Lett., {\bf 117} (2016), 113--117.

\bibitem{AD}
T.W.~Anderson, and D.A.~Darling, {\em Asymptotic theory of certain ``goodness of fit'' criteria based on stochastic processes}, Ann. Math. Stat., {\bf 23} (1952), N2, 193--212.


\bibitem{BNO}
L.~Beghin, Ya.~Nikitin, E.~Orsingher, 
{\em Exact small ball constants for some Gaussian processes under the 
$L_2$-norm}, ZNS POMI, {\bf298} (2003), 5--21;  J. Math. Sci., {\bf128} (2005), N1, 2493--2502.
 

\bibitem{DS}
N.~Dunford, J.T.~Schwartz, 
{\em Linear operators. Part III: Spectral operators}. 
With the assistance of William G. Bade and
Robert G. Bartle. Wiley-Interscience, N.Y., 1971.

\bibitem{F}
V.R. Fatalov, 
{\em Constants in the asymptotics of small deviation probabilities for {Gaussian} processes and fields},  Uspekhi mat. nauk, {\bf 58} (2003), N4(352), 89--134 (Russian); English transl.:
Russ. Math. Surv., {\bf58} (2003), 725--772.

\bibitem{GL}
F.~Gao, W.V.~Li, 
{\em Small ball probabilities for the Slepian Gaussian fields},
Trans. AMS, {\bf359} (2007), 1339--1350. 

\bibitem{Li}
W.V.~Li,
{\em Comparison results for the lower tail of Gaussian seminorms},
J.~Theor.\ Probab., {\bf5} (1992), N1, 1--31.

\bibitem{LS}
W.V.~Li, Q.M.~Shao,
{\em Gaussian processes: inequalities, small ball probabilities and 
applications}, Handbook Statist., {\bf19} (2001), 533--597.


\bibitem{Lf}
M.A.~Lifshits,
{\em Asymptotic behavior of small ball probabilities},
Prob.\ Theor.\ Math.\ Stat. (1999), Proc.\ VII Int.\ Vilnius Conference, 
453--468.

\bibitem{Lif2012-book}
M.A.~Lifshits, {\em Lectures on Gaussian processes}. Springer, 2012.

\bibitem{L-table}
M.A. Lifshits, 
{\em Bibliography of small deviation probabilities}, https://airtable.com/shrMG0nNxl9SiGxII/tbl7Xj1mZW2VuYurm, (2023).

\bibitem{Nm}
M.A.~Naimark, 
{\em Linear Differential Operators}. Ed.2. Moscow,
Nauka, 1969 (Russian); English transl. of the 1st ed.: Naimark
M.A. Linear Differential Operators. Part I (1967). 
N.Y.: F. Ungar Publ. Co. XIII. Part II (1968). 
N.Y.: F. Ungar Publ. Co. XV.

\bibitem{Na}
A.I.~Nazarov, 
{\em On the sharp constant in the small ball asymptotics of some 
Gaussian processes under $L_2$-norm}, Probl. Mat. Anal., {\bf26} (2003), 
179--214 (Russian); English transl.: J. Math. Sci., {\bf117} (2003), N3,
4185--4210. 

\bibitem{Naz09a}
A.I.~Nazarov, {\em Exact $L_2$-Small Ball Asymptotics of Gaussian Processes and the Spectrum of Boundary-Value Problems}, J. Theor. Probab., {\bf 22} (2009), N3, 640--665.


\bibitem{NN}
A.I.~Nazarov, Ya.Yu.~Nikitin,
{\em Exact $L_2$-small ball behavior of integrated
Gaussian processes and spectral asymptotics of boundary value problems},
Prob. Theor. Rel. Fields, {\bf129} (2004), N4, 469-494.

\bibitem{NN21}
A.I. Nazarov, Ya.Yu. Nikitin, {\em  Gaussian processes 
centered at their online average, and applications}, Stat. 
Prob. Lett., {\bf 170} (2021), 1--5. 

\bibitem{L2-survey}
A.I. Nazarov, Y.P. Petrova,
{\em  $L_2$-small ball asymptotics for Gaussian random functions: A survey}, Probab. Surv., {\bf 20} (2023), 608--663.

\bibitem{NP}
A.I. Nazarov, R.S. Pusev, {\em  Exact small deviation asymptotics in $L_2$-norm for some weighted Gaussian processes}, ZNS 
POMI, {\bf 364} (2009), 166--199 (Russian); 
English transl.: J. Math. Sci., {\bf 163} (2009), N4, 
409--429.

\bibitem{Sm} 
V.I.~Smirnov, 
{\em A course of higher mathematics}. Vol.IV, Part 2. Ed.6. Moscow, Nauka, 
1981 (Russian); English transl. of the 2nd ed.: Int. Series of 
Monographs in Pure and Appl. Math., Vol. 61, Oxford -- London -- 
New York -- Paris: Pergamon Press (1964).

\bibitem{T}
E.C.~Titchmarsh, 
{\em The theory of functions}, 2nd ed. London:
Oxford Univ. Press, 1975.

\bibitem{W}
G.S.~Watson, 
{\em Goodness-of-fit tests on a circle}, Biometrika, {\bf48} (1961), 109--114.


\end{thebibliography}
\end{document}